\documentclass{amsart}
\usepackage{mathrsfs,amsthm,amssymb,amsmath,url}

\theoremstyle{plain}
    \newtheorem{thm}{Theorem}
    
    \newtheorem{lem}[thm]{Lemma}
    \newtheorem{prop}[thm]{Proposition}
    \newtheorem{cor}[thm]{Corollary}
    
\theoremstyle{definition}
    \newtheorem{defn}[thm]{Definition}
    
\theoremstyle{remark}
    
\theoremstyle{remark}

\newcommand{\To}{\rightarrow}

\newcommand{\mult}{\times}
\newcommand{\ow}{\text{otherwise}}
\newcommand{\nin}{\notin}
\newcommand{\sm}{\setminus}

\newcommand{\cl}[1]{\langle #1 \rangle}

\renewcommand{\O}{{\mathscr O}}
\newcommand{\On}{{\mathscr O}^{(n)}}

\newcommand{\Oo}{{\mathscr O}^{(1)}}
\newcommand{\Ot}{{\mathscr O}^{(2)}}

\DeclareMathOperator{\pol}{Pol}

\newcommand{\inv}{^{-1}}

\newcommand{\C}{{\mathscr C}}

\newcommand{\F}{{\mathscr F}}

\newcommand{\I}{{\mathscr I}}
\newcommand{\A}{{\mathscr A}}
\newcommand{\B}{{\mathscr B}}

\newcommand{\E}{{\mathscr E}}
\newcommand{\G}{{\mathscr G}}
\renewcommand{\S}{{\mathscr S}}
\newcommand{\M}{{\mathscr M}}

\newcommand{\X}{{\mathscr X}}
\renewcommand{\H}{{\mathscr H}}

\newcommand{\un}{^{(n)}}

\newcommand{\uo}{^{(1)}}
\newcommand{\ut}{^{(2)}}

\newcommand{\rest}{\upharpoonright}

 \DeclareMathOperator{\id}{id}

\author[M.\,PINSKER]{MICHAEL PINSKER}
\address{Algebra\\TU Wien\\Wiedner Hauptstra\ss e 8-10/104\\A-1040 Wien, Austria}
\email{marula@gmx.at} \urladdr{http://www.dmg.tuwien.ac.at}
\title[Precomplete Clones Closed Under Conjugation]{Precomplete Clones On Infinite Sets Which Are Closed Under Conjugation}


\thanks{The author would like to thank L. Heindorf for his remarks on an earlier version of the paper.
Support through DOC [Doctoral Scholarship Programme of the Austrian
Academy of Sciences], and later through the Postdoctoral
Fellowship of the Japan Society for the Promotion of Science
(JSPS) is gratefully acknowledged}

\begin{document}

\begin{abstract}
    We show that on an infinite set,
    there exist no other precomplete clones closed under conjugation except those which
    contain all permutations. Since on base sets of some infinite cardinalities,
    in particular on countably infinite ones,
    the precomplete clones containing the permutations have been determined, this yields a complete list of the
    precomplete conjugation-closed clones in those cases. In addition, we show that there
    exist no precomplete submonoids of the full transformation monoid
    which are closed under conjugation except those which contain
    the permutations; the monoids of the latter kind are
    known.
\end{abstract}

\maketitle

\begin{section}{Background and the result}

    Let $X$ be a set and denote for all $n \geq 1$ the set of
    $n$-ary operations on $X$ by $\On$. The union $\O=\bigcup_{n\geq
    1}\On$ is the set of all operations on $X$ of finite arity. A
    \emph{clone} is a subset of $\O$ which contains all projections, i.e. all functions
    of the form $\pi^n_k(x_1,\ldots,x_n)=x_k$ ($1\leq k\leq n$), and
    which is closed under composition of functions. Ordering the
    clones on $X$ by set-theoretical inclusion, one obtains a
    complete algebraic lattice $Cl(X)$. We are interested in the
    structure of this
    lattice for infinite $X$, in which case it has cardinality
    $2^{2^{|X|}}$.

    We call a clone \emph{precomplete} or \emph{maximal} iff it is a dual atom in
    $Cl(X)$. The number of precomplete clones on an infinite base set
    equals the size of the whole clone lattice (\cite{Ros76}), and there is
    little hope to determine all of them. However, the precomplete
    clones which contain $\Oo$ have been determined on some infinite $X$ (\cite{Gav65}, \cite{GS02}).
    \begin{thm}\label{THM:maxaboveO1}
        If $X$ is countably infinite or of weakly compact
        cardinality, then there are exactly two precomplete clones $\pol(T_1)$ and $\pol(T_2)$
        above $\Oo$.
    \end{thm}

    For most other cardinalities of $X$, the number of precomplete
    clones above $\Oo$ is $2^{2^{|X|}}$, so in those cases it
    seems impossible to find them all (\cite{GS02}).

    The precomplete clones which contain the set of permutations $\S$ of the base set but not $\Oo$ have
    been determined on countably infinite $X$ in \cite{Hei02}, and we
    extended this result to all $X$ of regular cardinality in \cite{Pin033}. To
    describe these clones, the following concept was used:
    For a submonoid $\G\subseteq\Oo$, define the clone of
    \emph{polymorphisms} $\pol(\G)$ to consist of all $f\in\O$
    satisfying $f(g_1,\ldots,g_{n})\in\G$ whenever
    $g_1,\ldots,g_{n}\in\G$, where $n\geq 1$ is the arity of $f$. Call a subset $S\subseteq X$ \emph{large}
    iff $|S|=|X|$, and \emph{small} otherwise; $S$ is \emph{co-large (co-small)} iff its complement is large (small).
    A property holds for \emph{almost all} $x\in X$ iff there is a co-small $S\subseteq X$ such that the property holds for
    all $x\in S$. A function
    $f\in\Oo$ is \emph{almost surjective} iff almost all $x\in X$ are in the range of $f$. Then we have

    \begin{thm}\label{THM:bijections}
        Let $X$ be a set of regular cardinality $\kappa$. The precomplete clones
        over $X$ which contain all bijections but not all unary
        functions are exactly those of the form $\pol(\G)$, where
        $\G\in\{\A,\B,\E,\F\}\cup\{\G_\lambda: 1\leq\lambda\leq\kappa,\,\lambda\text{ a cardinal}\}$ is one of the following submonoids of $\Oo$:
        \begin{enumerate}
            \item{$\A=\{f\in\Oo:f^{-1}[\{y\}]$ is small for
            almost all $y\in X\}$}
            \item{$\B=\{f\in\Oo:f^{-1}[\{y\}]$ is small for all $y\in X\}$}
            \item{$\E=\{f\in\Oo:f$ is almost surjective$\}$}
            \item{$\F=\{f\in\Oo:f$ is almost surjective or constant$\}$}
            \item{$\G_\lambda=\{f\in\Oo:$ if $A\subseteq X$ has cardinality $\lambda$ then $|X\setminus f[X\setminus A]|\geq \lambda\}$}
        \end{enumerate}
    \end{thm}

    Clones containing the permutations $\S$ have the property that
    they are closed under conjugation, that is, $\C=\{\gamma\inv (f(\gamma(x_1),\ldots,\gamma(x_n))):
    f\in\C,\gamma\in\S\}$.
    We call clones with this property \emph{symmetric}; they are interesting
    because they are independent of the order or, indeed, any other structure that one might associate with
    the base set. But a clone need not contain
    $\S$ in order to be symmetric: For example, the clone consisting only of the projections and the constant functions
    is conjugation-closed.

    The set of symmetric clones is a sublattice $Cl_{sym}(X)$ of the
    clone lattice. Whereas $Cl(X)$ need not be dually atomic \cite{GS04}, that is, not every nontrivial clone
    is contained in a precomplete one, the sublattice of symmetric clones is. This is because
    there exist finitely many functions such that the only symmetric clone
    containing those functions is $\O$:
    If $\alpha\in\S$ is a permutation of $X$ which has large \emph{support}
    (i.e., if $\{x\in X:\alpha(x)\neq x\}$ is large),
    then $\alpha$ together
    with its conjugates (that is, all functions of the form $\gamma\inv\circ\alpha\circ\gamma$, where $\gamma\in\S$) generate $\S$ (\cite {SU33}
    for the countable and \cite{Bae34} for the uncountable case).
    And it is well-known that $\O$ is finitely generated over $\S$
    (see for example \cite{Hei02} for countably infinite and \cite{Pin033} for arbitrary infinite $X$).

    In the light of Theorems \ref{THM:maxaboveO1} and \ref{THM:bijections}, it is natural to ask whether
    it is possible to obtain a list of all symmetric precomplete
    clones not containing $\S$. We give the answer to this in
    this article.

    \begin{thm}\label{THM:allPrecompleteSymmClones}
        Let $X$ be infinite. If $\C$ is a symmetric precomplete clone, then it contains all permutations.
    \end{thm}
    We emphasize that this theorem is about symmetric clones
    which are dual atoms in $Cl(X)$, and not about symmetric clones which are dual atoms in the
    sublattice $Cl_{sym}(X)$ of symmetric clones. There exist clones of
    the latter type which do not contain all permutations: For
    example, the clone consisting of all $f\in\O$ for which the
    set
    $\{x\in X: f(x,\ldots,x)=x\}$ is co-small does not contain $\S$, is obviously symmetric,
    and it follows easily from
    the complete description of the clone lattice above this clone in \cite{GS0y} that
    there is no non-trivial symmetric clone containing it.

    \begin{cor}
        If $X$ has regular cardinality, then the symmetric precomplete
        clones which do not contain $\Oo$ are exactly those from Theorem
        \ref{THM:bijections}. If $X$ is
        countably infinite or of weakly compact cardinality, then the symmetric
        precomplete clones are exactly those from Theorems
        \ref{THM:maxaboveO1} and
        \ref{THM:bijections}.
    \end{cor}

    A \emph{unary clone} is a clone consisting only of \emph{essentially unary} functions, where by an essentially
    unary function we mean one which depends on only one of its variables. Clearly, though formally different,
    unary clones can be seen as submonoids of the full transformation monoid
    $\Oo$ and we shall not distinguish between the two notions.
    A \emph{precomplete} unary clone is a dual atom in the lattice of submonoids of $\Oo$. In
    \cite{Gav65}, all precomplete unary
    clones which contain all permutations were determined on countably infinite $X$. The result was
    generalized in
    \cite{Pin033} to arbitrary infinite sets: For a cardinal $1\leq\lambda\leq |X|$ we call a function
    $f\in\Oo$ \emph{$\lambda$-injective} iff there exists
    $A\subseteq X$ with $|A|<\lambda$ such that the
    restriction of $f$ to $X\sm A$ is injective.

    \begin{thm}\label{THM:allPrecompleteMonoids}
        Let $X$ be a set of infinite cardinality $\kappa$. If $\kappa$ is regular, then the
        precomplete submonoids of $\Oo$ which contain the permutations
        are exactly the monoid $\A$ and the monoids $\G_\lambda$ and $\M_\lambda$ for
        $\lambda=1$ and $\aleph_0\leq \lambda\leq\kappa$, $\lambda$ a cardinal, where
        \begin{itemize}

            \item{$\A=\{f\in\Oo:f^{-1}[\{y\}]$ is small for
            almost all $y\in X\}$}
            \item{$\G_\lambda=\{f\in\Oo:$ if $A\subseteq X$ has cardinality $\lambda$
            then $|X\setminus f[X\setminus A]|\geq \lambda\}$}
            \item{$\M_\lambda=\{f\in\Oo: f \text{ is
            }\lambda\text{-surjective or not
            }\lambda\text{-injective}\}$}
        \end{itemize}

        If $\kappa$ is singular, then the same is true
        with the monoid $\A$ replaced by
        $$
            \A'=\{f\in\Oo: \exists \lambda < \kappa \,\,(\,|f\inv[\{x\}]|\leq\lambda\text{ for almost all }x\in X\,)\,\}.
        $$
    \end{thm}

    We shall obtain the following result.

    \begin{thm}\label{THM:allPrecompleteSymmMonoids}
        Let $X$ be infinite. The symmetric precomplete unary clones are exactly those of Theorem
        \ref{THM:allPrecompleteMonoids}.
    \end{thm}

    It might be interesting to note that on finite $X$, all
    symmetric clones are known (\cite{Kho92},\cite{Kho93},\cite{Kho94},\cite{Mar96a},\cite{Mar96b},
    see also the survey paper \cite{Sze03}). If $X$ has at least five elements, then
    the only symmetric precomplete
    clone is the S\l upecki clone of all
    functions which are either essentially unary or take at most $|X|-1$
    values. In this case the
    clone of all idempotent functions is the only other
    clone which is maximal in $Cl_{sym}(X)$ (but not in $Cl(X)$, since it is properly contained in
    the clone of all $f\in\O$ for which $f(a,\ldots,a)=a$, for any fixed $a\in X$). For $|X|<5$ the situation is more complicated, see
    \cite{Sze03}. On all finite $X$ with at least three elements, there exist precisely two
    symmetric precomplete submonoids of $\Oo$: The first one contains all permutations,
    and is the monoid of all functions $f\in\Oo$ which are either a
    permutation or for which the set $\{x\in X:
    \exists y\neq x: f(x)=f(y)\}$ has at least three elements. The
    second one does not contain the permutations and consists of the non-permutations plus all even
    permutations. That there are no other symmetric precomplete monoids, even no other dual atoms in the lattice
    of symmetric monoids, is not difficult
    to prove from the definitions
    of these two monoids and the fact (\cite{LM91}) that if
    $\M$ is a symmetric monoid on a finite set $X$, then
    $\M\cup\S$ is a monoid which contains $\M$.

    Since on an infinite set $X$ of size $\kappa=\aleph_\alpha$ there
    exist $\max\{2^{2^{|\alpha|}},2^{2^{\aleph_0}}\}$ submonoids
    of $\Oo$ containing $\S$ (\cite{Pin041}), describing all symmetric monoids or even clones seems to
    be hopeless.
    \subsection{Notation}
    For a set of functions $\H\subseteq\O$, we use
    the common notation $\cl{\H}$ for the smallest clone containing
    $\H$. If $n\geq 1$, then $\H\un$ denotes the $n$-ary functions
    in $\H$. In particular, if $\H\subseteq\Oo$ is a set of unary functions, then $\cl{\H}\uo$ is the unary part
    of the clone generated by $\H$ and therefore nothing else than the monoid generated by $\H$.\\
    The unary projection $\pi^1_1$ is the identity function on $X$
    and we denote it also by $\id$. For $S\subseteq
    X$ we denote the image of $S$ under $f$ by $f[S]$.
\end{section}

%----------------------------------------------------------------------------------------------------------

\begin{section}{Unary clones: The proof of Theorem \ref{THM:allPrecompleteSymmMonoids}}

    This section contains the proof of the simple unary case (Theorem \ref{THM:allPrecompleteSymmMonoids}).
    We start by citing a theorem which is essential for the whole paper; for us,
    its most important implication is that a symmetric clone which
    does not contain all permutations in fact contains almost
    no permutations.
    \begin{thm}[\cite{SU33},\cite{Bae34}]\label{THM:SchreierUlam}
        If $\alpha\in\S$ has large support, then $\alpha$ together
        with its conjugates generate $\S$.
    \end{thm}
    \begin{cor}\label{COR:noPermWithLargeSupport}
        If $\C$ is a symmetric clone which does not contain
        all permutations, then $\C$ does not contain any
        permutation with large support.
    \end{cor}
    \begin{lem}\label{LEM:representationUnary}
        If $\G$ is a symmetric submonoid of $\Oo$ and
        $\alpha\in\S$, then $\cl{\{\alpha\}\cup\G}\uo=\{\alpha^n\circ
        g: g\in\G, n\geq 0\}=\{g\circ\alpha^n: g\in\G, n\geq 0\}$.
    \end{lem}
    \begin{proof}
        We begin with the first equality. It is clear that $\{\alpha^n\circ
        g: g\in\G, n\geq 0\}\subseteq \cl{\{\alpha\}\cup\G}\uo$. The
        other inclusion we prove by induction over terms $t$ in $\cl{\{\alpha\}\cup\G}\uo$.
        The statement is obvious for $t=\alpha$ and $t=g\in\G$. So
        assume that $t=\alpha\circ s$, where $s\in\cl{\{\alpha\}\cup\G}\uo$ satisfies the
        induction hypothesis. Then there exist $n\geq 0$ and $g\in
        \G$ with $s=\alpha^n\circ g$ so that we have $t=\alpha^{n+1}\circ
        g$. To finish the induction, assume that $t=h\circ s$,
        with $h\in \G$ and $s\in\cl{\{\alpha\}\cup\G}\uo$ satisfying the induction hypothesis, so that $s=\alpha^n\circ g$ for some $n\geq 0$
        and some $g\in\G$. Set $h'=\alpha^{-n}\circ
        h\circ\alpha^n$. Then $h'\in\G$ because $\G$ is symmetric,
        and $h=\alpha^n\circ h'\circ \alpha^{-n}$. Hence,
        $t=h\circ\alpha^n\circ s=\alpha^n\circ h'\circ
        \alpha^{-n}\circ \alpha^n\circ g=\alpha^n\circ (h'\circ
        g)$ and we are finished.\\
        Now for the second equality, observe that $t=\alpha^n\circ g$, with $n\geq 0$ and $g\in \G$,
        if and only if
        $t=g'\circ\alpha^n$, where $g'=\alpha^{n}\circ g\circ \alpha^{-n}\in\G$.
    \end{proof}

    \begin{proof}[Proof of Theorem \ref{THM:allPrecompleteSymmMonoids}]
        We show that if $\G$ is a symmetric submonoid of $\Oo$, and if
        $\G\nsupseteq\S$, then $\G$ is not
        precomplete.
        Take $\alpha\in\S$ with large and co-large support, and take
        $\beta\in\S$ with large and co-small support.
        Having large support, neither $\alpha$ nor $\beta$ are elements of
        $\G$, by Corollary \ref{COR:noPermWithLargeSupport}. If $\G$ was a precomplete monoid,
        then every function in $\Oo$ would be an element of $\cl{\{\alpha\}\cup\G}\uo$, which equals
        $\{\alpha^n\circ
        g: g\in\G, n\geq 0\}$ by Lemma \ref{LEM:representationUnary}. In particular, there would
        exist $n\geq 1$ and $\gamma\in\G$ such that
        $\beta=\alpha^n\circ\gamma$. Obviously, $\gamma$ has to be
        a permutation; being an element of $\G$, by Corollary \ref{COR:noPermWithLargeSupport} it must have
        small support. The power $\alpha^n$ still has co-large
        support, and so does $\alpha^n\circ\gamma$. Hence,
        $\beta$ cannot equal $\alpha^n\circ\gamma$, so that we end up with a contradiction.
    \end{proof}
\end{section}

%---------------------------------------------------------------------------------------------------------------

\begin{section}{Non-unary clones: The proof of Theorem \ref{THM:allPrecompleteSymmClones}}

    We turn to the proof of Theorem
    \ref{THM:allPrecompleteSymmClones}. The reason why this is more
    difficult is that with non-unary functions, Lemma
    \ref{LEM:representationUnary} becomes more complicated.

    \begin{lem}\label{LEM:representationOfArbitrary}
        If $\C$ is a symmetric clone and
        $\alpha\in\S$ and $t\in\On$, then $t\in
        \cl{\{\alpha\}\cup\C}$ iff there exists some integer $k\geq 1$, $f\in\C^{(n\cdot k)}$,
        and a sequence $(a_{i,j}:1\leq i\leq n,1\leq j\leq k)$ of nonnegative integers such that
        $$
            t(x_1,\ldots,x_n)=f(\alpha^{a_{1,1}}(x_1),\ldots,\alpha^{a_{1,k}}(x_1),\ldots,\alpha^{a_{n,1}}(x_n),
            \ldots,\alpha^{a_{n,k}}(x_n)).
        $$
    \end{lem}
    \begin{proof}
        It is clear that if $t$ is of that form, then $t\in\cl{\{\alpha\}\cup\C}$. The
        other direction we prove by induction over terms $t$ in $\cl{\{\alpha\}\cup\C}$.
        The statement is obvious for $t=\alpha$ and $t=f\in\C$. So
        assume that $t=\alpha\circ s$, where $s\in\cl{\{\alpha\}\cup\C}$ satisfies the
        induction hypothesis. Then $s=f(\alpha^{a_{1,1}}(x_1),\ldots,\alpha^{a_{1,k}}(x_1),\ldots,
        \alpha^{a_{n,1}}(x_n),\ldots,\alpha^{a_{n,k}}(x_n))$
        for some $f\in\C$, $k\geq 1$ and a sequence $(a_{i,j})$ of integers.
        Set $f'=\alpha\circ f\circ\alpha\inv\in \C$. Then
        $f=\alpha\inv\circ f'\circ\alpha$ and we can calculate
        $$
        \begin{aligned}
            t&=\alpha\circ s\\
            &=\alpha\circ\alpha\inv\circ
            f'\circ\alpha(\alpha^{a_{1,1}}(x_1),\ldots,\alpha^{a_{1,k}}(x_1),\ldots,\alpha^{a_{n,1}}(x_n),\ldots,
        \alpha^{a_{n,k}}(x_n))\\
            &=f'(\alpha^{a_{1,1}+1}(x_1),\ldots,\alpha^{a_{1,k}+1}(x_1),\ldots,\alpha^{a_{n,1}+1}(x_n),\ldots,
            \alpha^{a_{n,k}+1}(x_n)).
        \end{aligned}
        $$
        To finish the induction, assume that $t=f(t_1,\ldots,t_{m})$,
        where $f\in \C^{(m)}$ and $t_l\in\cl{\{\alpha\}\cup\C}^{(n)}$
        satisfies the induction hypothesis, $1\leq l\leq
        m$. Then
        $$
            t_l=f_l(\alpha^{a^l_{1,1}}(x_1),
            \ldots,\alpha^{a^l_{1,k_l}}(x_1),\ldots,\alpha^{a^l_{n,1}}(x_n),
            \ldots,\alpha^{a^l_{n,k_l}}(x_n))
        $$
        for $f_l\in\C$, $k_l\geq 1$ and sequences
        $(a^l_{i,j}:1\leq i\leq n,1\leq j\leq k_l)$. By adding fictitious variables
        to the $f_l$, which we can do freely within
        a clone, we can assume that all $k_l$ are equal to some $k\geq 1$, and
        even that all sequences $(a^l_{i,j}:1\leq i\leq n,1\leq j\leq
        k_l)$ are identical with one sequence $(a_{i,j}:1\leq i\leq n,1\leq j\leq
        k)$. Then the $f_l$ are all of the same arity
        $k\cdot n$. Setting $g=f(f_1,\ldots,f_m)\in\C$, we obtain
        $$
            t=g(\alpha^{a_{1,1}}(x_1),
            \ldots,\alpha^{a_{1,k}}(x_1),\ldots,\alpha^{a_{n,1}}(x_n),
            \ldots,\alpha^{a_{n,k}}(x_n)).
        $$
    \end{proof}

    To make things more convenient, we consider permutations
    $\alpha$ which satisfy $\alpha^2=\id$; then the preceding
    lemma becomes
    \begin{lem}\label{LEM:reprArbitraryWithNiceAlpha}
        If $\C$ is a symmetric clone and
        $\alpha\in\S$ satisfying $\alpha^2=\id$, and if $t\in\On$, then $t\in
        \cl{\{\alpha\}\cup\C}$ iff there exists $f\in\C^{(2n)}$
        such that
        $$
            t(x_1,\ldots,x_n)=f(\alpha
            (x_1),\ldots,\alpha(x_n),x_1,\ldots,x_n).
        $$
    \end{lem}
    \begin{proof}
        If $t$ has such a representation, then it obviously is an element $\cl{\{\alpha\}\cup\C}$. Conversely, let
        $t\in \cl{\{\alpha\}\cup\C}$. By the preceding lemma, $t$
        has a representation
        $$
            t(x_1,\ldots,x_n)=f(\alpha^{a_{1,1}}(x_1),\ldots,\alpha^{a_{1,k}}(x_1),\ldots,\alpha^{a_{n,1}}(x_n),
            \ldots,\alpha^{a_{n,k}}(x_n))
        $$
        for some $k\geq 1$, a sequence $(a_{i,j}:1\leq i\leq n,1\leq
        j\leq k$), and $f\in\C$.
        Now because $\alpha^2=\id$, all the exponents $a_{i,j}$ become
        either $0$ or $1$. Since a clone is closed under
        identification of variables and under changing of the
        order of variables, as well as under addition of fictitious variables, we can assume
        that $\alpha^1(x_i)$ and $\alpha^0(x_i)$ occur exactly
        once in the representation for each $1\leq i\leq n$, and
        that they occur in the desired order.
    \end{proof}

    \begin{defn}
        Let $n\geq 1$. We say that $\H\subseteq\O$ has the \emph{$n$-ary co-large approximation
        property} iff for all $f\in\On$ and all co-large $S\subseteq X$
        there exists $g\in\H\un$ such that
        $g\rest_{S^{n}}=f\rest_{S^{n}}$. We call the
        function $g$ an \emph{approximation} to $f$ on $S$. $\H$ has the
        \emph{co-large approximation property} iff it has the $n$-ary
        co-large approximation property for all $n\geq 1$.
    \end{defn}

    A reformulation of the definition of the co-large approximation property which we will use heavily in the
    following is: $\H\subseteq\O$ has the $n$-ary co-large approximation property iff for all co-large $S\subseteq X$
    and all functions $f: S^n\To X$ there is a function $g\in\H$ extending $f$ to $X^n$. We will give an example of a non-trivial clone having the
    co-large approximation property later in this paper (Proposition \ref{PROP:PolEandPolF}). The
    following is not surprising since all finitary operations are generated by binary operations.

    \begin{lem}\label{LEM:binaryApproxProperty}
        Let $\C\subseteq\O$ be a clone. Then $\C$ has the binary co-large
        approximation property iff it has the co-large approximation property.
    \end{lem}
    \begin{proof}
        It suffices to show that the binary co-large approximation property implies the
        co-large approximation property. Let $f\in\On$ and $S\subseteq X$ co-large. We want to find in $\C$ an
        approximation to $f$ on $S$. Without loss of generality, we can assume $S$ to be large. Fix a bijection $j:X\To S$ and set $h=j\circ
        f$. Then $h[X^n]\subseteq S$,
        hence $h\rest_{S^n}$ is an operation on $S$. Because on every set the binary
        functions generate all finitary functions, $h\rest_{S^n}$ is generated by binary functions on
        $S$, all of which have extensions to $X$ in
        $\C$ by the binary co-large approximation property. Hence,
        there exists $g\in\C\un$ such that $g\rest_{S^n}=h\rest_{S^n}$.
        Moreover, again by the binary co-large approximation property of $\C$,
        there is $t\in\C\uo$ such that $t\rest_S=j\inv$.
        Thus for $t\circ g\in\C$ we have $t\circ g\rest_{S^n}=j\inv\circ j\circ
        f\rest_{S^n}=f\rest_{S^n}$, so $t\circ g$ is an approximation to $f$ on $S$.
    \end{proof}

    \begin{lem}\label{LEM:approxArbitraryFunctions}
        Let $\C$ be symmetric and precomplete, and assume that $\C\nsupseteq\S$. Then $\C$ has the
        co-large approximation property.
    \end{lem}
    \begin{proof}
        Let $f\in\On$ and any co-large $S\subseteq X$ be given. Take
        $\alpha\in\S$ with support $X\sm S$ and such that $\alpha^2=\id$, i.e., let $\alpha$ be the identity
        on $S$ and let its cycles on $X\sm S$ be any partition of $X\sm S$ into two-element sets.
        Because $\alpha$ has large support, Corollary \ref{COR:noPermWithLargeSupport} implies that
        $\alpha\nin\C$,
        and therefore $\cl{\{\alpha\}\cup\C}=\O$ as $\C$ is precomplete. By Lemma \ref{LEM:reprArbitraryWithNiceAlpha},
        there
        exists $g\in\C^{(2n)}$ such that
        $f(x_1,\ldots,x_n)=g(\alpha(x_1),\ldots,\alpha(x_n),x_1,\ldots,x_n)$.
        Set $h(x_1,\ldots,x_n)=g(x_1,\ldots,x_n,x_1,\ldots,x_n)\in\C$. Since
        $\alpha$ is the identity on $S$, we have $h(x_1,\ldots,x_n)=g(x_1,\ldots,x_n,x_1,\ldots,x_n)=
        g(\alpha(x_1),\ldots,\alpha(x_n),x_1,\ldots,x_n)=f(x_1,\ldots,x_n)$ on $S$, whence $h$ is an approximation
        to $f$ on $S$.
    \end{proof}

    \begin{defn}
        We say that a function $f\in\Oo$ is \emph{generous} iff
        all equivalence classes of its kernel are large. We set
        $\I_0$ to consist of all generous functions which are
        onto.
    \end{defn}

    The unary clone $\I_0$ is an example of a symmetric proper
    submonoid of $\Oo$ having the unary co-large approximation
    property. Indeed, given any $f\in\Oo$ and any co-large $S\subseteq X$,
    then since $X\sm S$ is large we can find $\gamma: X\sm S \To X$ mapping $X\sm S$ onto $X$
    in such a way that $\gamma\inv [\{y\}]$ is large for all
    $y\in X$. Now set $g(x)=f(x)$ if $x\in S$, and
    $g(x)=\gamma(x)$ if $x\nin S$. Then $g\inv[\{y\}]$ is large for all $y\in X$ since $g$ extends $\gamma$,
    so $g\in\I_0$; moreover, $g$ is an approximation to $f$ on $S$
    by its definition.\\

    Let $f\in\I_0$, and fix any $x\in X$. Then we have the following
    possibilities:
    \begin{itemize}
        \item{There exist $n\geq 0$ and $p\geq 1$ such that $f^n(x)=f^{n+p}(x)$.}
        \item{$x,f(x),f^2(x),\ldots$ are all distinct.}
    \end{itemize}

    In the first case, choosing $p$ minimal with the property that there is $n\geq 0$ such that
    $f^{n}(x)=f^{n+p}(x)$,
    we call the set of all elements $y\in X$ such that there is $k\geq 0$ with
    $f^k(y)=x$ or $f^k(x)=y$, or equivalently the connectedness component of $x$ in
    the graph $(X,f)$, a \emph{$p$-snail}.
    In the second case, we call the set of all elements of $X$ connected to $x$ a
    \emph{$0$-snail}. This definition is independent of the
    element of a snail we are looking at, that is, it depends only on the connectedness component of
    $x$ in $(X,f)$.
    Indeed, let $y\in X$ be connected to $x$; then
    there exists $k\geq 0$ such that either $f^k(y)=x$
    or $f^k(x)=y$. Assume without loss of generality
    the first case holds. Now if $f^{n}(x)=f^{n+p}(x)$
    for some $n\geq 0$ and $p\geq 1$, and if $p$ is minimal with this
    property, then $f^{k+n}(y)=f^{k+n+p}(y)$,
    so in particular $y$ does not think it is part of an infinite snail. Suppose there
    exist $j\geq 0$ and $1\leq q<p$ such that
    $f^j(y)=f^{j+q}(y)$. Then
    $f^{j+k}(y)=f^{j+q+k}(y)$ and so
    $f^{j}(x)=f^{j+q}(x)$, contradicting the minimality
    of $p$. Hence, $y$ thinks it is in a $p$-snail too. If on the other hand $x,f(x),f^2(x),\ldots$
    are all distinct then the same holds for $y$, for if $f^j(y)=f^{j+q}(y)$ for $j\geq 0$ and $q\geq 1$,
    then $f^j(x)=f^{j+k}(y)=f^{j+q+k}(y)=f^{j+q}(x)$,
    contradiction.

    \begin{defn}
        We call a function $f\in\I_0$ \emph{rich} iff it has a large number of $p$-snails for all
        $p\geq 0$.
    \end{defn}

    If $(a_p)_{p\in\omega}$ is any sequence of cardinals
    $\leq\kappa$
    which is not constantly zero, then there exists a function
    $f\in\I_0$ whose number of $p$-snails is $a_p$, for all
    $p\in\omega$. Indeed, let $p>0$, take
    $x_0,\ldots,x_{p-1}\in X$, and define $f(x_i)=x_{i+1}$ for
    $0\leq i\leq p-2$, and $f(x_{p-1})=x_0$. Moreover, let $f$ map
    $X\sm\{x_0,\ldots,x_{p-1}\}$ onto $X$ in such a way that the
    preimage of every $y\in X$ is large. Now set $X'$ to consist
    of all $x\in X$ such that there is $k\geq 0$ with
    $f^k(x)\in\{x_0,\ldots,x_{p-1}\}$. Then the restriction of $f$ to $X'$ is a
    function on $X'$ which is onto and generous, and which has exactly one $p$-snail and no other
    snails. Since $|X'|=|X|$, such a function exists also on $X$.
    Similarly we can prove the existence of functions having
    only one $0$-snail and no other snails. By taking the union
    over functions on disjoint sets, all of which have only one snail, according to the sequence
    $(a_p)_{p\in\omega}$, one obtains the function whose number of $p$-snails is $a_p$,
    for all $p\in\omega$. In
    particular, rich functions exist.

    \begin{lem}\label{LEM:existsRich}
        Let $\H\subseteq\O$ have the unary co-large approximation property. Then
        there exists a rich $f\in\H\uo$.
    \end{lem}
    \begin{proof}
        Fix any large and co-large $S\subseteq X$, and some $T\subseteq S$
        such that $T$ and $S\sm T$ are large. Take any $m_1:
        T\To T$ which is rich (as an operation on the base set
        $T$). Now let $m_2: S\sm T\To X$ be so that $m_2\inv
        [\{y\}]$ is large for every $y\in X$. Set $m=m_1\cup m_2: S\To X$; by
        the unary co-large approximation property of $\H$, there is a function
        $f\in\H\uo$ with $f\rest_S=m$. Because $f$ extends $m_2$, every $y\in X$ has a large preimage under $f$;
        thus, $f\in\I_0$. Being identical with $m_1$ on $T$, $f$ has a large
        number of $p$-snails for all $p\geq 0$, and we see that $f$ is rich.
    \end{proof}

    \begin{lem}\label{LEM:isomorphicGraphs}
        Let $f,g\in\I_0$. Then $f$ and $g$ have the same number of
        $p$-snails for all $p\geq 0$ iff there exists $\gamma\in\S$ such
        that $f=\gamma\inv\circ g\circ\gamma$.
    \end{lem}
    \begin{proof}
        If $f=\gamma\inv\circ g\circ\gamma$, then the structures
        $(X,f)$ and $(X,g)$ are isomorphic via $\gamma$; this obviously implies
        that if $x\in X$ is part of a $p$-snail of $f$, then
        $\gamma(x)$ is part of a $p$-snail of $g$. In particular,
        $f$ and $g$ have the same number of $p$-snails for all
        $p\geq 0$.\\
        Let on the other hand $f$ and $g$
        have the same number of $p$-snails for all $p\geq 0$; we will construct an isomorphism
        $\gamma: (X,f)\To (X,g)$. Assume first that $f$ and $g$
        have exactly one $p$-snail and no other snails, where $p>0$. There exist
        $a,b\in X$ such that $f^p(a)=a$ and $g^p(b)=b$. Set $A=\{a,f(a),\ldots,f^{p-1}(a)\}$ and
        $B=\{b,g(b),\ldots,g^{p-1}(b)\}$; then, since we are in a $p$-snail, $|A|=|B|=p$. Now for every
        $x\in X$ there is a minimal $n\geq 0$ such that
        $f^n(x)\in A$; we say that $x$ is on the $n$-th level with respect to $f$. Levels
        with respect to $g$ are defined analogously. We
        define $\gamma$ by induction over levels. Set
        $\gamma(f^k(a))=g^k(b)$, for all $0\leq k\leq p-1$; that defines $\gamma$ on $A$, which constitutes
        exactly the $0$-th
        level of $f$. To define $\gamma$ on the
        first level, consider the sets $f\inv[\{x\}]\sm A$ and $g\inv[\{\gamma(x)\}]\sm B$,
        for all $x\in A$. Since
        they are both large, we can map
        the first bijectively onto the latter; we extend $\gamma$
        by such a bijection. This defines $\gamma$ for level $1$. Say $\gamma$ has
        already been defined for all $x$ of a level smaller than
        $n$, where $n\geq 2$. Let $x$ be at level $n-1$ and take again the sets $f\inv[\{x\}]$ and
        $g\inv[\{\gamma(x)\}]$; this time we do not have to remove the elements of $A$ and $B$,
        since they are not
        mapped to level $n-1$, as $n\geq 2$. The preimages are both large, so as on the first level, we map
        the first bijectively onto the latter, and extend $\gamma$
        by such a bijection. This defines $\gamma$ for level $n$, since every $y$ on level $n$ is an element of
        the preimage of exactly one $x$ of level $n-1$.
        Moreover, because every $x\in X$ appears on exactly one level, $\gamma$ is a function on
        $X$, and it maps the $n$-th level with respect to $f$
        to the $n$-th level with respect to $g$, for all $n\geq 0$. We claim that $\gamma$ is injective.
        Indeed, to
        see this we prove by induction over levels that $\gamma$ is injective on each level, which is sufficient as
        it maps distinct levels of $f$ to distinct levels of $g$.
        It is clear that $\gamma$ is injective on $A$, by definition. So let $x,y$ be distinct elements of the
        $n$-th level, where $n\geq 1$. If
        $f(x)\neq f(y)$, then by induction hypothesis
        $\gamma(f(x))\neq \gamma(f(y))$ and so
        $\gamma(x)\neq\gamma(y)$ since $\gamma(x)\in
        g\inv[\{\gamma(f(x))\}]$ but $\gamma(y)\in
        g\inv[\{\gamma(f(y))\}]$. If $f(x)=f(y)$, then, since both
        $x$ and $y$ are in $f\inv[\{f(x)\}]\sm A$, and since
        $\gamma$ maps $f\inv[\{f(x)\}]\sm A$ bijectively onto
        $g\inv[\{\gamma(f(x))\}]\sm B$, we again have that
        $\gamma(x)\neq\gamma(y)$, and the induction is complete. To prove that $\gamma$ is
        surjective, we again proceed by induction. If $y\in B$,
        then it is in the range of $\gamma$. Now let $y$ be on
        level $n$ with respect to $g$, where $n\geq 1$,
        and assume that all elements of lower level are in the range of $\gamma$;
        then $g(y)\in\gamma[X]$ by
        induction hypothesis. Let $z$ be the element for which $\gamma(z)=g(y)$. By definition of $g$,
        $f\inv[\{z\}]\sm A$ is mapped onto
        $g\inv[\{\gamma(z)\}]\sm B=g\inv[\{g(y)\}]\sm B$; therefore, $y$ has
        an element in $f\inv[\{z\}]$ mapped to it, so $g$ is surjective.
        Now $\gamma(f(x))=g(\gamma(x))$ for all $x\in X$ by
        construction of $\gamma$, so $\gamma$ is an isomorphism.\\
        Assume next that $f$ and $g$ have only one $0$-snail, and no other
        snails. Take any $a,b\in X$, and set
        $A=\{a,f(a),f^2(a),\ldots\}$ and
        $B=\{b,g(b),g^2(b),\ldots\}$. Define
        $\gamma(f^k(a))=g^k(b)$, for all $k\geq 0$; then proceed
        again defining $\gamma$ by induction over levels. As before,
        one checks
        that $\gamma$ is an isomorphism.\\
        Now if $f$ and $g$ have the same number of $p$-snails for
        all $p\geq 0$, then fix for
        every $p\geq 0$ a bijection $\sigma_p$ from the $p$-snails of $f$
        onto the $p$-snails of $g$. By the preceding discussion,
        if $\partial$ is any $p$-snail of $f$, we can construct an
        isomorphism $\gamma_\partial$ from $(\partial,f)$ onto
        $(\sigma_p(\partial),g)$. If we set $\gamma$ to be the
        union over all $\gamma_\partial$, for all snails
        $\partial$ of $f$, we obtain an isomorphism between
        $(X,f)$ and $(X,g)$.
    \end{proof}

    \begin{lem}\label{LEM:allRich}
        If $\C$ is a symmetric clone with the unary co-large approximation property,
        then $\C$ contains all rich functions.
    \end{lem}
    \begin{proof}
        By Lemma \ref{LEM:existsRich}, $\C$ contains a rich
        function. Hence it contains all rich functions
        by the preceding lemma.
    \end{proof}

    Recall from Theorem \ref{THM:bijections} that $\E\subseteq\Oo$ consists
    of those functions in $\Oo$ which take all but a small set of values.
    \begin{lem}\label{LEM:underPolE}
        If $\C$ is a symmetric clone which has the unary co-large approximation property,
        and if $\C\uo\subseteq\E$, then
        $\C\subseteq\pol(\E)$.
    \end{lem}
    \begin{proof}
        Assuming that there exists $f\in\C\sm\pol(\E)$ we find $g\in\C\uo$ such that $g\nin\E$.
        Because $f\nin\pol(\E)$, there
        exist $\alpha_1,\ldots,\alpha_n\in\E$ such that
        $f(\alpha_1,\ldots,\alpha_n)\nin\E$. Fix $A\subseteq X$ small such
        that $A\cup\alpha_i[X]=X$ for all $1\leq i\leq n$; this is possible since the range of all
        $\alpha_i$ is co-small. Choose
        $S\subseteq X\sm A$ large and co-large such that
        $f(\alpha_1,\ldots,\alpha_n)[X]\cup\{f(x,\ldots,x):x\in S\}$ is
        still co-large. Set $Y=X\sm(A\cup S)$. Then $Y$ is large, implying that we can find $\gamma\in\Oo$
        mapping $Y$ onto $X$. Now define for all $1\leq i\leq n$
        $$
            \beta_i(x)=\begin{cases}\alpha_i\circ\gamma(x),&x\in
            Y\\x,&x\in A\cup S.\end{cases}
        $$
        Then
        $f(\beta_1,\ldots,\beta_n)[X]=f(\beta_1,\ldots,\beta_n)[Y]\cup
        f(\beta_1,\ldots,\beta_n)[A\cup S]=f(\alpha_1,\ldots,\alpha_n)[X]\cup\{f(x,\ldots,x):x\in A\cup
        S\}$ is co-large. Let $h\in\Oo$be so that $h[S]=S$, $h\rest_S$ is
        a rich function on the large set $S$, and
        such that it maps $X\sm S$ onto $X$ in such a way that every $x\in X$
        has a large preimage in $X\sm S$ under $h$.
        Then already
        its definition on $X\sm S$ guarantees that $h$ is onto and that
        all classes of its kernel are large, so $h\in\I_0$.
        Because $h$ is onto we have that for all $1\leq i\leq n$,
        $\beta_i\circ h[X]=\beta_i[X]=\beta_i[Y]\cup A\cup S=\alpha_i[X]\cup A\cup S=X$. Also, since all classes
        in the kernel of $h$ are large, so are those of $\beta_i\circ h$;
        therefore, $\beta_i\circ h\in\I_0$. Moreover, $\beta_i\circ
        h$ has a large number of $p$-snails for all $p\geq 0$,
        since already $h\rest_S$ has this property and
        $\beta_i\circ h\rest_S=h\rest_S$ as $\beta_i\rest_S=\id_S$.
        Hence, $\beta_i\circ h$ is rich and therefore an element of $\C$ by Lemma \ref{LEM:allRich}.
        But $f(\beta_1\circ h,\ldots,\beta_n\circ h)[X]= f(\beta_1,\ldots,\beta_n)[X]$ is
        co-large so that it suffices to set $g=f(\beta_1\circ h,\ldots,\beta_n\circ
        h)\in\C$.
    \end{proof}
    Remember that the monoid $\F\subseteq\Oo$ is the union of $\E$ and
    all constant functions.
    \begin{lem}\label{LEM:underPolF}
        If $\C$ is a symmetric clone which has the unary co-large approximation property,
        and if $\C\uo\nsubseteq\E$ and $\C\uo\subseteq\F$, then
        $\C\subseteq\pol(\F)$.
    \end{lem}
    \begin{proof}
        A slight modification of the proof of the preceding
        lemma yields that assuming there is $f\in\C\sm\pol(\F)$, we can find $g\in\C\uo$ such that
        $g\nin\F$. So let $f\in\C\sm\pol(\F)$; then there exist
        $\alpha_1,\ldots,\alpha_n\in\F$ such that
        $f(\alpha_1,\ldots,\alpha_n)\nin\F$. Now observe that the
        conditions $\C\uo\nsubseteq\E$ and $\C\uo\subseteq\F$ imply that $\C\uo$
        contains a constant function; hence it contains all
        constant functions as it is symmetric. The functions
        $\alpha_1,\ldots,\alpha_n\in\F$ are either constant or
        almost surjective; assume $\alpha_1,\ldots,\alpha_k$ are constant, and $\alpha_{k+1},\ldots,\alpha_n$
        almost surjective, where $1\leq k< n$. Note that $k=n$ is impossible for otherwise $f(\alpha_1,\ldots,\alpha_n)$
        would be constant and thereby an element of $\F$. Consider $f(\alpha_1,\ldots,\alpha_k,x_{k+1},\ldots,x_n)$ as
        a $(n-k)$-ary function $\tilde{f}$ of the variables $x_{k+1},\ldots,x_n$. By the proof of
        the preceding lemma, we can find rich functions
        $\beta_{i}\circ h$ for $k+1\leq i\leq n$ such that
        $g=\tilde{f}(\beta_{k+1}\circ h,\ldots,\beta_{n}\circ
        h)$ has co-large range. It also follows from that proof that
        $g[X]\supseteq
        \tilde{f}(\alpha_{k+1},\ldots,\alpha_{n})[X]=
        f(\alpha_{1},\ldots,\alpha_{n})[X]$, and so $g$ is not constant. Therefore $g\nin \F$. Now
        it is enough to observe that
        $g=f(\alpha_1,\ldots,\alpha_k,\beta_{k+1}\circ h,\ldots,\beta_{n}\circ
        h)$, and that all functions which appear here as arguments of
        $f$ are either constant or rich and thus elements of
        $\C$, by Lemma \ref{LEM:allRich}. Whence, $g\in\C\uo$,
        contradicting the assumption $\C\uo\subseteq\F$.
    \end{proof}
    \begin{lem}\label{LEM:containsChi}
        If $\C$ is a symmetric clone which has the unary co-large approximation property,
        and if $\C\uo\nsubseteq\F$, then
        $\X=\{\rho\in\Oo:|\rho[X]|\leq 2\text{ and }\rho \text{ is
         generous}\}$ is contained in $\C$.
    \end{lem}
    \begin{proof}
        To start with, let $\rho\in\X$ be so that it takes two distinct values $c_1,c_2\in X$,
        and write $\theta_i=\rho\inv[\{c_i\}]$, $i=1,2$.
        Take any $f\in\C\uo\sm\F$. Since $f$ is not constant there exist
        $a_1\neq a_2$ in the range of $f$. Define $s: f[X]\To
        \{a_1,a_2\}$ by $s(a_1)=a_1$ and $s(x)=a_2$ for all $x\neq a_1$. Since $f\nin\F$ we have that $f[X]$ is
        co-large, and so
        the unary co-large approximation property of $\C$ implies that we can
        extend $s$ to a function $g\in\C\uo$. Choose any rich function
        $h\in\Oo$; then $h\in\C$ by Lemma \ref{LEM:allRich}.
        Therefore, $g\circ f\circ
        h$ is an element of $\C$ as well. Now since $h$ is generous, so is $g\circ f\circ h$.
        Also, since $h$ is onto and by the construction of $g$, $g\circ f\circ h[X]=g\circ
        f[X]=\{a_1,a_2\}$. Write $\zeta_i=(g\circ f\circ h)\inv[\{a_i\}]$, $i=1,2$. Since the $\zeta_i$ and the $\theta_i$
        are large, there exists
        $\gamma\in\S$ mapping $\theta_i$ onto $\zeta_i$, for $i=1,2$.
        Then $\gamma\inv\circ g\circ f\circ h\circ\gamma$ is in $\C$ as $\C$ is symmetric, and maps
        $\theta_i$ to $\gamma\inv(a_i)$, $i=1,2$. Now let
        $t\in\C\uo$ be so that it maps $\gamma\inv(a_i)$ to $c_i$,
        $i=1,2$. We can find such a $t$ in $\C$ by the unary co-large
        approximation property. Then $t\circ\gamma\inv\circ g\circ f\circ h\circ\gamma\in
        \C$, and it maps $\theta_i$ to $c_i$, $i=1,2$. Whence, it
        equals $\rho$ so that we infer $\rho\in\C$. Therefore $\C$ contains all functions in $\X$ which
        take two values.\\
        Now to see that $\C$ contains the constant functions as well, let $c\in X$, and let $f$ be as before.
        By the unary co-large approximation property, we can find $q\in\C\uo$ mapping all
        elements of the co-large
        range of $f$ to $c$. Then $q\circ f\in \C$ is constant with value $c$ and we are done.
    \end{proof}

    \begin{lem}\label{LEM:containsLargeAndColarge}
        If $\C$ is a symmetric clone which has the co-large approximation property, and if $\C\supseteq\X$, then there exists $g\in\C\uo$
        having
        large and co-large range.
    \end{lem}
    \begin{proof}
        Choose any large and co-large $T\subseteq X$, and any element $0\in T$.
        Fix $f\in\C\ut$ such that $f[\{0\}\mult X]$
        and $f[X \mult \{0\}]$ are large, and such that $f\rest_{T^2}$ is constantly
        $0$. To obtain $f$, let $S\subseteq X\sm T$ be large and so
        that $T\cup S$ is co-large, and let $\gamma: S\To X$ be onto.
        Then define a partial binary operation $m$ to yield
        constantly $0$ on $T^2$, and to satisfy $m(s,0)=m(0,s)=\gamma(s)$ for all $s\in S$. Since the domain
        of $m$ is contained in $(S\cup T)^2$ and $S\cup T$ is co-large, by the
        co-large approximation property we can extend $m$ to
        $f\in\C\ut$. Clearly, $f$ yields constantly $0$ on $T^2$ as $m$
        does, and $f[\{0\}\mult X]\supseteq f[\{0\}\mult S]=m[\{0\}\mult
        S]=\gamma[S]=X$, and the same holds for $f[X\mult \{0\}]$,
        so $f$ is indeed as desired.
        We distinguish two cases.\\
        \textbf{Case 1.} For all $c\in X$ it
        is true that $f[X\mult\{c\}]$ and $f[\{c\}\mult X]$ are
        co-small. Then consider an arbitrary large and co-large  $A\subseteq
        X$ with $0\nin A$. Set
        $\Gamma=f\inv[X\sm A]\subseteq X^2$ and let
        $\alpha: X\sm T\To\Gamma$ be onto. By the assumption for this case,
        $f[X\mult\{c\}]\sm A$ and $f[\{c\}\mult X]\sm A$ are still
        large for all $c\in X$. Thus
        the components $\alpha_i=\pi^2_i\circ\alpha$ hit every $c\in X$ at a large number of arguments,
        $i=1,2$. Moreover, by taking the union with any rich function on the base set $T$,
        we can extend the $\alpha_i$ to $T$ so that $\alpha_i[T]=T$ and
        $\alpha_i\rest_T$ is rich. The resulting operations $\alpha_i\in\Oo$
        still hit every $c\in X$ at a large number of
        arguments, so they are elements of $\I_0$. Also, since they already have a large number of $p$-snails on
        $T$ for all $p\geq 0$, they certainly still have this property
        on $X$. Therefore they are rich
        and hence elements of $\C\uo$ by Lemma \ref{LEM:allRich}. But now setting
        $g=f(\alpha_1,\alpha_2)\in\C$ we have that $g[X]=g[X\sm T]\cup
        g[T]=f(\alpha_1,\alpha_2)[X\sm T]\cup \{0\}=f[X^2]\sm A$ is large and
        co-large and we are done.\\
        \textbf{Case 2.} There exists $c\in X$ such that either $f[X\mult\{c\}]$ or $f[\{c\}\mult X]$
        is co-large,
        say without loss of generality this is the case for $f[\{c\}\mult X]$. Since $f[\{0\}\mult X]$
        is large we can choose
        $\Gamma\subseteq X$ large and co-large such that
        $f[\{0\}\mult\Gamma]$ is large and such that $f[\{c\}\mult
        X]\cup f[\{0\}\mult\Gamma]$ is still
        co-large. Take moreover a rich $\beta\in\I_0$; then $\beta\in\C$ by
        Lemma \ref{LEM:allRich}. Now we define
        $\alpha\in\Oo$ by
        $$
            \alpha(x)=
            \begin{cases}0,&\beta(x)\in\Gamma\\
            c,&\ow.
            \end{cases}
        $$
        The range of $\alpha$ equals $\{0,c\}$, and the preimage of both values under $\alpha$ is large.
        Hence $\alpha\in\X\subseteq\C$. Thus it is enough to set
        $g=f(\alpha,\beta)\in\C$ and observe
        that $g[X]=f[\{c\}\mult (X\sm\Gamma)]\cup f[\{0\}\mult\Gamma]$
        is large and co-large.

    \end{proof}

    \begin{lem}\label{LEM:largeAndColargeEverything}
        Assume that $\C$ is a symmetric clone with the co-large approximation property. If
        there is $g\in\C\uo$ with large and co-large range, then
        $\C=\O$.
    \end{lem}
    \begin{proof}
        Because $\C$ has the co-large approximation property, by
        Lemma \ref{LEM:existsRich}
        there exists a rich $h\in\C\uo$. The function $g'=g\circ
        h\in\C\uo$ is generous and has large and co-large range.
        Set $S=g'[X]\subseteq X$.
        There exists $\gamma\in\S$ with the property that
        for all distinct $x,y\in X$ it is true that $g'(x)=g'(y)$ implies
        $g'\circ\gamma(x)\neq g'\circ\gamma(y)$. Indeed, let $\{\theta_i\}_{i\in\kappa}$ be an enumeration of
        the classes of the
        kernel of $g'$, and enumerate the elements of those classes by
        $\theta_i=\{x_i^j\}_{j\in\kappa}$, for all $i\in\kappa$. Now define $\gamma$ by $\gamma(x_i^j)=x_j^i$ for all
        $i,j\in\kappa$. Every $x\in X$ is equal to some $x_j^i$, and therefore has the element $x_i^j$ mapped
        to it by $\gamma$; hence $\gamma$ is surjective. If $x_i^j\neq x_p^q$, then $i\neq p$ or $j\neq q$
        and thus
        $\gamma(x_i^j)=x_j^i\neq x_q^p=\gamma(x_p^q)$, and we see that $\gamma$ is injective.
        Now if $g'(x)=g'(y)$ for
        $x\neq y$, then
        $x,y\in\theta_i$ for some $i\in\kappa$, and so there are distinct $j,k\in\kappa$ such that $x=x_i^j$ and $y=x_i^k$.
        But then $g'(\gamma(x))=g'(x_j^i)$ and $g'(\gamma(y))=g'(x_k^i)$ are not equal as $x_j^i$ and
        $x_k^i$ belong to
        different kernel classes of $g'$. Thus $\gamma$ has the desired properties. Set
        $g''=\gamma\inv\circ g'\circ\gamma$; then $g''\in\C$ because $\C$ is symmetric, and
        $g''$ still satisfies $g''(x)\neq g''(y)$ whenever $g'(x)=g'(y)$ and
        $x,y\in X$ are distinct. Since the range
        of $g''$ is large and co-large, the co-large approximation property guarantees
        that we can find $f\in\C\uo$ which
        maps $g''[X]$ injectively onto $S$. Set $g'''=f\circ g''$.
        Then the function $t(x)=(g'(x),g'''(x))$ maps $X$
        injectively into $S^2$. For both $g'$ and $g'''$ take only values in $S$, and if $g'(x)=g'(y)$
        for distinct
        $x,y\in X$, then $g'''(x)\neq g'''(y)$. Therefore, the function $s(x,y)=(g'(x),g'''(x),g'(y),g'''(y))$
        maps $X^2$ injectively into $S^4$. Now let an arbitrary $q\in\O\ut$ be
        given. We can find
        $m\in\C^{(4)}$ satisfying
        $$
            m(g'(x),g'''(x),g'(y),g'''(y))=q(x,y)
        $$
        for all $x,y\in X$. This is because distinct pairs in $X^2$ yield distinct quadruples in $S^4$
        via $s$, and we can define $m$ on $S^4$ as required by $q$ and extend it to a function in $\C^{(4)}$
        by the co-large approximation property. The equality implies $q\in\C$
        and so
        $\C\supseteq\O\ut$ since $q$ was arbitrary. Whence,
        $\C=\O$.
    \end{proof}
    \begin{prop}\label{PROP:PolEandPolF}
        $\pol(\E)$ and $\pol(\F)$ have the co-large approximation
        property.
    \end{prop}
    \begin{proof}
        Let $f\in\Ot$ and $S\subseteq X$
        be co-large. We construct an approximation $g$ to $f$ on $S$
        which is an element of both $\pol(\E)$ and $\pol(\F)$,
        proving that these clones have the binary co-large approximation
        property and hence the co-large approximation property by
        Lemma
        \ref{LEM:binaryApproxProperty}.
        Let $\gamma\in\Oo$ be so that it maps $X\sm S$ onto $X$. Define
        $g\in\Ot$ by
        $$
            g(x_1,x_2)=
            \begin{cases}
                f(x_1,x_2),&x_1,x_2\in S\\
                \gamma(x_1),&x_1\nin S\\
                \gamma(x_2),&x_1\in S\wedge x_2\nin S.\\
            \end{cases}
        $$
        Then $g$ and $f$ agree on $S^2$ so that $g$ is an approximation to $f$ on $S$. To see that
        $g\in\pol(\E)$, let $g_1,g_2\in\E$ be given. Since the
        range of $g_1$ is co-small, we have that
        $g_1$ misses only a small number of elements of $X\sm S$, and therefore
        $\gamma[(X\sm S)\cap
        g_1[X]]$ is co-small. Now
        $g(g_1,g_2)(x)=\gamma(g_1(x))$ whenever $g_1(x)\in X\sm
        S$. Hence, $g(g_1,g_2)[X]\supseteq \gamma[(X\sm S)\cap
        g_1[X]]$ is co-small, so that $g(g_1,g_2)\in\E$. Whence,
        $g\in\pol(\E)$ and $\pol(\E)$ has the co-large approximation property.\\
        We now show that $g\in\pol(\F)$. Let $g_1,g_2\in\F$, i.e. they are either constant or almost surjective.
        If both $g_1$ and $g_2$ are almost surjective, then so is
        $g(g_1,g_2)$ as we have just seen, so $g(g_1,g_2)\in\F$.
        If on the other hand both functions are constant, then the
        composite $g(g_1,g_2)$ is constant as well so that again
        $g(g_1,g_2)\in\F$. Next assume that $g_2$ is constant and $g_1$ is almost
        surjective. Then exactly the argument of the preceding
        paragraph yields that $g(g_1,g_2)$ is almost surjective,
        since in that argument we needed only that $g_1$ is almost
        surjective, and no assumptions on $g_2$. Therefore we
        obtain $g(g_1,g_2)\in\F$ also in this case. Finally,
        consider the case where $g_1$ is constant and $g_2$ is
        almost surjective. We distinguish two subcases: If $g_1$
        constantly yields a value $c$ in $X\sm S$, then
        $g(g_1,g_2)(x)=\gamma(c)$ for all $x\in X$, so
        $g(g_1,g_2)$ is constant. If on the other hand $g_1$ is constant with value
        $c\in S$, then by definition of $g$ we have that
        $g(g_1,g_2)(x)=\gamma(g_2(x))$ whenever $g_2(x)\in X\sm
        S$. Since $g_2$ is almost surjective, it misses only a
        small number of elements in $X\sm S$, and so $\gamma[(X\sm S)\cap
        g_2[X]]$ is co-small. Hence, since $g(g_1,g_2)[X]\supseteq\gamma[(X\sm S)\cap
        g_2[X]]$, we see that $g(g_1,g_2)$ is almost surjective. In either
        case, $g(g_1,g_2)\in\F$ and thus $g\in\pol(\F)$. Therefore
        $\pol(\F)$ has the co-large approximation property.
    \end{proof}

    \begin{prop}\label{PROP:onlyPolEandPolF}
        The only symmetric precomplete clones having the
        co-large approximation property are $\pol(\E)$ and $\pol(\F)$.
    \end{prop}
    \begin{proof}
        We know from Theorem \ref{THM:bijections} that $\pol(\E)$ and $\pol(\F)$ are precomplete and
        symmetric; by the preceding proposition, both clones have the co-large approximation property.
        Suppose $\C$ is a symmetric clone having the co-large approximation property, and which is distinct from those two clones. If $\C\uo\subseteq\E$, then
        $\C\subseteq\pol(\E)$ by Lemma \ref{LEM:underPolE} and so $\C$ is not precomplete.
        Moreover, if
        $\C\uo\subseteq\F$ and $\C\nsubseteq\E$, then $\C\subseteq\pol(\F)$ by Lemma
        \ref{LEM:underPolF}, hence $\C$ is not precomplete either.
        Finally, if $\C\uo\nsubseteq\F$, then $\C\supseteq\X$ by Lemma \ref{LEM:containsChi}, and thus
        $\C$ contains a unary function with large and co-large
        range by Lemma \ref{LEM:containsLargeAndColarge}. Hence,
        $\C=\O$ by Lemma \ref{LEM:largeAndColargeEverything}.
    \end{proof}

    \begin{proof}[Proof of Theorem \ref{THM:allPrecompleteSymmClones}]
        Assume there exists a symmetric and precomplete clone $\C$ with
        $\C\nsupseteq\S$. Then $\C$ has the co-large approximation property by Lemma
        \ref{LEM:approxArbitraryFunctions}, so $\C$ has to equal
        either $\pol(\E)$ or $\pol(\F)$ by the preceding proposition. Contradiction.
    \end{proof}
\end{section}

\end{document}